\documentclass[reqno,draft]{amsart}

\usepackage{graphicx}
\usepackage[nobysame]{amsrefs}
\usepackage{hyperref}

\newcommand{\R}{{\mathbb R}}

\newcommand{\F}{{\mathcal F}}
\newcommand{\G}{{\mathcal G}}
\newcommand{\T}{{\mathcal T}}
\newcommand{\al}{\alpha}
\newcommand{\be}{\beta}

\newcommand{\eps}{\varepsilon}

\newtheorem{prop}{Proposition}
\newtheorem{theorem}{Theorem}
\newtheorem{corollary}{Corollary}
\newtheorem*{fht}{Fractional Helly Theorem}
\newtheorem*{fhtb}{Fractional Helly Theorem for boxes}

\theoremstyle{definition}
\newtheorem{example}{Example}
\newtheorem{case}{Case}

\title[A fractional Helly theorem for boxes]{A fractional Helly theorem for boxes}

\author{I. B\'ar\'any}
\address{Alfr\'ed R\'enyi Institute of Mathematics,
PO Box 127, H-1364 Budapest, Hungary,
and Department of Mathematics, University College London,
Gower Street, London, WC1E 6BT, U.K.}
\email{barany@renyi.hu}

\author{F. Fodor}
\address{Department of Geometry, Bolyai Institute, University of
Szeged, Aradi v\'ertan\'uk tere 1, H-6720 Szeged, Hungary,
and Department of Mathematics and Statistics,
University of Calgary, Canada}
\email{fodorf@math.u-szeged.hu}

\author{A. Mart\'{\i}nez-P\'{e}rez}
\address{Universidad de Castilla- La Mancha
Departamento de An\'{a}lisis Econ\'{o}mico y Finanzas.
Universidad de Castilla-La Mancha. 
Avda. Real F\'{a}brica de Seda, s/n. 45600
Talavera de la Reina. Toledo.  Spain.}
\email{alvaro.martinezperez@uclm.es}

\author{L. Montejano}
\address{Instituto de Matem\'{a}ticas, Universidad Nacional Aut\'onoma de
M\'exico, \'Area de la Investigaci\'on Cient\'\i fica, Circuito Exterior, Cu.
Coyoacan 04510, M\'exico D.F., M\'exico}
\email{luis@matem.unam.mx}

\author{D. Oliveros}
\address{Instituto de Matem\'{a}ticas, Universidad Nacional Aut\'onoma de
M\'exico, \'Area de la Investigaci\'on Cient\'\i fica, Circuito Exterior, Cu.
Coyoacan 04510, M\'exico D.F., M\'exico}
\email{dolivero@matem.unam.mx}

\author{A. P\'or}
\address{Department of Mathematics, Western Kentucky University, Bowling Green,
KY 42101, USA}
\email{attila.por@wku.edu}

\dedicatory{This paper is dedicated to Javier Bracho on occasion of his sixtieth birthday.}

\begin{document}
\begin{abstract}
Let $\F$ be a family of $n$ axis-parallel boxes in $\R^d$ 
and $\al\in (1-1/d,1]$ a real number. There exists a real number $\be(\al)>0$ such that 
if there are $\al{n\choose 2}$ intersecting
pairs in $\F$, then $\F$ contains an intersecting
subfamily of size $\be n$. A simple example shows that 
the above statement is best possible in the sense
that if $\al\leq 1-1/d$, then
there may be no point in $\R^d$ that belongs to more than $d$ elements of $\F$. 
\end{abstract}

\maketitle

\section{Introduction and results}
According to the classical theorem of Helly \cite{DGK}, if every $d+1$-element subfamily
of a finite family of convex sets in $\mathbb{R}^{d}$ has nonempty intersection, 
then the entire family has nonempty intersection. Although the 
number $d+1$ in Helly's theorem cannot be lowered in general, 
it can be reduced for some special families of convex sets. For example, 
if any two elements in a finite family of axis-parallel boxes in $\R^d$ 
intersect, then all members of the family intersect, cf. \cite{DG}. 

Katchalski and Liu \cite{KL1979} proved the following generalization 
of Helly's theorem for the case when not all 
but only a fraction of $d+1$-element subfamilies have a nonempty intersection in a family
of convex sets.

\begin{fht}\label{th:fracH}{\rm (Katchalski and Liu~\cite{KL1979})}
Assume that $\al \in (0,1]$ is a real number and $\F$ is a family of $n$ convex
sets in $\R^d$.
If at least $\al\binom n{d+1}$ of the $(d+1)$-tuples of $\F$ intersect,
then $\F$ contains an intersecting subfamily of size $\frac {\al}{d+1} n$.
\end{fht}

The bound on the size of the intersecting subfamily was later improved by Kalai \cite{K}
from $\frac {\al}{d+1} n$ to $(1-(1-\al)^{1/(d+1)})n$, and this bound is best possible. 

In this paper, we study the fractional behaviour of finite 
families of axis-parallel boxes, or boxes for short.
We note that the boxes can be either open or closed, our statements
hold for both cases. 
Our aim is to prove a statement similar to the 
Fractional Helly Theorem.

The intersection graph $\G_{\F}$ of a finite family $\F$ of boxes
is a graph whose vertex set is the set of elements of $\F$, and
two vertices are connected by an edge in $\G_{\F}$ precisely when the 
corresponding boxes in $\F$ have nonempty intersection. 

Recall that for two integers $n\geq m\geq 1$, the Tur\'an-graph $\T(n,m)$ is a complete
$m$-partite graph on $n$ vertices in which the cardinalities of the 
$m$ vertex classes are as close to each other as possible.  Let 
$t(n,m)$ denote the number of edges of the Tur\'an graph $\T(n,m)$.
It is known that $t(n,m)\leq (1-\frac{1}{m})\frac{n^2}{2}$, and
equality holds if $m$ divides $n$. Furthermore,
\begin{equation}\label{turangraph}
\lim_{n\to\infty}\frac{t(n,m)}{\frac{n^2}{2}}=1-\frac{1}{m}. 
\end{equation}
For more information on the properties of Tur\'an graphs see, for example,
the book of Diestel \cite{Diestel2010}.

The following example shows that we cannot hope for a 
statement for boxes that is completely analogous to the 
Fractional Helly Theorem.

\begin{example}\label{example}
Let $n\geq d+1$ and $m, k\geq 0$ be integers such that
$n=md+k$ and $0\leq k\leq d-1$. Let $n_1,\ldots, n_d$ be positive integers with 
$n=n_1+\cdots +n_d$ and $n_i=\lceil\frac{n}{d}\rceil$ for $1\leq i\leq k$ and
$n_i=\lfloor\frac{n}{d}\rfloor$ for $k+1\leq i\leq d$.
 For $1\leq i\leq d$, consider $n_i-1$ 
hyperplanes orthogonal to the $i$th coordinate direction. 
These hyperplanes cut $\R^d$ into $n_i$ pairwise disjoint open slabs 
$B'_{ij}, j=1,\ldots, n_i$.  
Let $C$ be a large open axis-parallel box that 
intersects each slab and let $\F_i$ consist of the open boxes 
$B_{ij}=C\cap B'_{ij}$. 
Define $\F$ as the union of the $\F_i$. 
   
This way we have obtained a family $\F$ of $n$ boxes with the property 
that two elements of $\F$ intersect exactly if they belong to different $\F_i$.  
The intersection graph of $\F$ is $\T(n,d)$ and thus the number of 
intersecting pairs in $\F$ is $t(n,d)$. However, there is no point
of $\R^d$ that belongs to any $d+1$-element subfamily of $\F$. Thus,
\eqref{turangraph} shows that in a fractional Helly-type statement for 
boxes, the percentage $\al$ has to be greater than $1-\frac{1}{d}$.     
\end{example}

Let $n\geq k\geq d$ and let $T(n,k,d)$ denote the maximal
number of intersecting pairs in a family $\mathcal{F}$ of $n$ boxes in 
$\R^d$ with the property that no $k+1$ boxes in  $\F$ have a point in common.

\begin{theorem}\label{turanindex} 
With the above notation,
\begin{equation*}
T(n,k,d)<\frac{d-1}{2d}n^{2}+\frac{2k+d}{2d}n.
\end{equation*}
\end{theorem}

It is quite easy to precisely determine $T(n,k,d)$ when $d=1$:

\begin{prop}\label{intervals} $T(n,k,1)=(k-1)n-{k \choose 2}$.
\end{prop}

Theorem~\ref{turanindex} directly implies the following corollary. 

\begin{corollary}\label{thm:frac-boxes}
Assume that $\eps>0$ is a real number and $\F$ is 
a family of $n$  boxes in $\R^{d}$. If at least 
$\left (\frac{d-1}{2d}+\eps \right )n^{2}$ pairs of $\F$ intersect,
then $\F$ contains an intersecting subfamily of
size $dn\eps -\frac{d}{2}+1$.
\end{corollary}

The proof of Corollary~\ref{thm:frac-boxes} is given in Subsection~2.2.
Corollary~\ref{thm:frac-boxes} yields the next theorem, which is our main result.

\begin{fhtb}
For every $\al\in (1-\frac{1}{d},1]$ there exists a real number $\be(\al)>0$
such that, for every family $\F$ of $n$ boxes in $\R^d$, if an $\al$
fraction of pairs are intersecting in $\F$, 
then $\F$ has an intersecting subfamily of cardinality at least $\be n$. 
\end{fhtb}

Kalai's lower bound $\be(\al)=1-(1-\al)^{1/(d+1)}$ 
for the size of the intersecting subfamily in 
the fractional Helly theorem yields that if $\al\to 1$, then 
$\beta(\al)\to 1$ as well. The same holds for families of
parallel boxes as stated in the following theorem.

\begin{theorem}\label{thm:limit}
Let $\F$ be a family of $n$ boxes in $\R^d$, 
and let $\al \in (1-\frac{1}{d^2},1]$ be a real number. 
If at least $\al{n\choose 2}$ pairs of boxes in $\F$ intersect, 
then there exists a point that belongs to at least 
$(1-d\sqrt{1-\al})n$ elements of $\F$.
\end{theorem}

Simple calculations show that Corollary~\ref{thm:frac-boxes} does not
imply Theorem~\ref{thm:limit} so we provide a separate proof for it 
in Section~2.

\section{Proofs}
\subsection{Proof of Theorem~\ref{turanindex}}
It is enough to prove that if  no $k+1$ elements of $\mathcal{F}$ 
have a point in common, then there are at least $\frac{n^{2}-2(k+d)n}{2d}$ 
non-intersecting pairs. We may assume by standard
arguments that the boxes in $\mathcal{F}$ are all open, so $B\in \mathcal{F}$ is of the form 
$B=(a_{1}(B),b_{1}(B))\times \dots \times (a_{d}(B),b_{d}(B))$. We assume
without loss of generality that all numbers $a_{i}(B),b_{i}(B)$ ($B\in \mathcal{F}$)
are distinct. For $B\in \mathcal{F}$ we define $\deg B$ to be the number of boxes in 
$\mathcal{F} $ that intersect $B$.

We prove Theorem~\ref{turanindex} by induction on $n$. The starting case $n=k$ is simple
since then $\frac{n^2-2(k+d)n}{2d}<0$. In the induction step $n-1 \to n$ we
consider two cases.

\begin{case} {\em When there is a box $B$ with 
$\deg B\leq (1-\frac{1}{d})n+\frac{2k+1}{2d}$.} 

By induction, we have at least 
$\frac{(n-1)^{2}-2(k+d)(n-1)}{2d}$ non-intersecting pairs after removing $B$ from 
$\mathcal{F} $. Since $B$ is involved in at least $(n-1)-\left( 1-\frac{1}{d}\right) n-
\frac{2k+1}{2d}$ non-intersecting pairs, there are at least 
\[
\frac{(n-1)^{2}-2(k+d)(n-1)}{2d}-1+\frac{n}{d}-\frac{2k+1}{2d}=\frac{
n^{2}-2(k+d)n}{2d} 
\]
non-intersecting pairs in $\mathcal{F}$, indeed.
\end{case}

\begin{case} {\em For every 
$B\in \mathcal{F}$ $\deg B\geq (1-\frac{1}{d})n+\frac{2k+1}{2d}$.} 

We show by contradiction that this cannot happen which finishes the
proof.

We define $d$ distinct boxes $B_1,\ldots,B_d \in \mathcal{F}$ the following way. Set 
$$c_1=\min \{b_1(B): B \in \mathcal{F}\}$$ 
and define $B_1$ via $c_1=b_1(B_1)$. The box 
$B_1$ is uniquely determined as all $b_1(B)$ are distinct numbers. Assume now
that $i<d$ and that the numbers $c_1,\ldots,c_{i-1}$, and boxes 
$B_1,\ldots,B_{i-1}$ have been defined. Set 
$$c_i=\min \{b_i(B): B \in \mathcal{F}
\setminus\{B_1,\ldots,B_{i-1}\}\}$$ 
and define $B_i$ via $c_i=b_i(B_i)$ which
is unique, again.

Let $\mathcal{F}^{\prime }= \mathcal{F} \setminus \{ B_1, \dots, B_d\}$. We partition $\mathcal{F}^{\prime
}$ into $d+2$ parts. Let $\mathcal{F}_0$ be the set of all boxes of $\mathcal{F}^{\prime }$ that
intersect every $B_i$. For $i=1, \dots, d$ let $\mathcal{F}_i$ be the set of all boxes
in $\mathcal{F}^{\prime }$ that intersect every $B_j$ for $j \neq i$ but do not
intersect $B_i$. Let $\mathcal{F}^*$ be the set of all boxes of $\mathcal{F}^{\prime }$ that
intersect at most $d-2$ of the $B_i$ boxes. As this is a partition of $\mathcal{F}^{\prime }$ we have 
\[
|\mathcal{F}_0| + \sum_{i=1}^{d} |\mathcal{F}_i| + |\mathcal{F}^*| = |\mathcal{F}^{\prime }| = n-d. 
\]
Note that $|\mathcal{F}_0| \leq k$ since every box in $\mathcal{F}_0$ contains the point $(c_1,
\dots, c_d)$.

Let $N$ be the number of intersecting pairs between $\{B_1,\ldots,B_d\}$ and 
$\mathcal{F}^{\prime }$. Each $B_i$ intersects at least $\deg B_i -(d-1)$ boxes from $\mathcal{F}^{\prime }$ as 
$B_i$ may intersect $B_j$ for all $j\in [d], j\ne i$. Since
every $\deg B_i\ge (1-\frac{1}{d})n +\frac{2k+1}{2d}$ we have 
\[
d \left ((1-\frac{1}{d})n+\frac{2k+1}{2d} - (d-1)\right ) \le N 
\]

Every box in $\mathcal{F}_{0}$ intersects every $B_{i}$, $i\in \lbrack d]$, every box
in $\mathcal{F}_{i}$ intersects every $B_{j}$ except for $B_{i}$ and every box in 
$\mathcal{F}^{\ast }$ intersects at most $(d-2)$ of the $B_{i}$. Consequently 
\[
N\leq d|\mathcal{F}_{0}|+(d-1)\sum_{i=1}^{d}|\mathcal{F}_{i}|+(d-2)|\mathcal{F}^{\ast }|. 
\]%
So we have 
\begin{eqnarray*}
d\left( (1-\frac{1}{d})n+\frac{2k+1}{2d}-(d-1)\right) &\leq
&d|\mathcal{F}_{0}|+(d-1)\sum_{i=1}^{d}|\mathcal{F}_{i}|+(d-2)|\mathcal{F}^{\ast }| \\
&=&|\mathcal{F}_{0}|+(d-1)\left( |\mathcal{F}_{0}|+\sum_{1}^{d}|\mathcal{F}_{i}|+|\mathcal{F}^{\ast }|\right)
-|\mathcal{F}^{\ast }| \\
&=&|\mathcal{F}_{0}|+(d-1)(n-d)-|\mathcal{F}^{\ast }|.
\end{eqnarray*}%
Simplifying the inequality and using $|\mathcal{F}_{0}|\leq k$ give 
\[
k+\frac{1}{2}\leq |\mathcal{F}_{0}|-|\mathcal{F}^{\ast }|\leq k-|\mathcal{F}^{\ast }| 
\]%
implying $|\mathcal{F}^{\ast }|\leq -\frac{1}{2}$, which is a contradiction.
\end{case}

\subsection{Proof of Corollary~\ref{thm:frac-boxes}}
If no point of $\R^d$ belongs to $dn\eps-\frac{d}{2}+1$ elements of $\F$, then 
by Theorem~\ref{turanindex} the number of intersecting pairs of $\mathcal{F}$ is smaller than 
$$\frac{d-1}{2d}n^{2}+\frac{2(dn\eps -\frac{d}{2})+d}{2d}n=\left (\frac{d-1}{2d}+\eps \right )n^{2},$$
which yields a contradiction.

\subsection{Proof of Theorem~\ref{thm:limit}}
Let $\pi _{i}$ denote the orthogonal projection to the $i$th dimension in $\R^d$, that is, 
$\pi _{i}(B)=(a_{i}(B),b_{i}(B))$ for $B\in\F$. Set $\varepsilon
=1-\alpha $. Define $T_{i}=\{\pi _{i}(B):B\in \F\}$; this is a family of $n$
intervals, and all but at most $\eps \binom{n}{2}$ of the pairs in 
$T_{i}$ intersect. According to the sharp version of the fractional
Helly theorem (cf. \cite{K}), $T_{i}$ contains an intersecting subfamily $T_{i}^{\prime
}$ of size $(1-\sqrt{\varepsilon })n$, let $c_{i}$ be a common point of all
the intervals in $T_{i}^{\prime }$. Define $D_{i}=\{B\in \F:c_{i}\notin \pi
_{i}(B)\}$. Then $\mathcal{F}\setminus \bigcup_{1}^{d}D_{i}$ consists of at least 
$(1-d\sqrt{\varepsilon })n=(1-d\sqrt{1-\alpha })n$ boxes and all of them contain
the point $(c_{1},\ldots ,c_{d})$.

\subsection{Proof of Proposition~\ref{intervals}}
Let $k\in\{1,\dots,n\}$ be an integer, and let $\F$ be the family of open intervals $(i,i+k)$ for $i=1,2,\ldots,n$. Thus $\F$ consists of $n$ intervals, no $k+1$ of them have a point in common, and there are $(k-1)n-{k \choose 2}$ intersecting pairs in $\F$. Consequently $T(n,k,1)\geq (k-1)n-{k \choose 2}$.

Next we show, by induction on $n$ that $T(n,k,1)\leq (k-1)n-{k \choose 2}$.
Let $\F$ be a family of $n$ intervals such that no $k+1$ of them have a common point. We assume that these intervals are closed which is no loss of generality. The statement is clearly true when $n=k$.
Let $[a,b] \in F$ be the interval where $b$ is minimal. 
Since any interval intersecting $[a,b]$ contains $b$, there are at most $k-1$ intervals intersecting
$[a,b]$. Removing $[a,b]$ from $\F$ and applying induction, we find there are at most $(k-1)(n-1)-{k \choose 2}$ intersecting pairs in $\F\setminus \{[a,b]\}$. That is, there are at most $k-1+(k-1)(n-1)-{k\choose 2} = (k-1)n-{k\choose 2}$ intersecting pairs in $\F$.

\section{Acknowledgements}
The authors wish to acknowledge the support of this research by the Hungarian-Mexican
Intergovernmental S\&T Cooperation Programme grant T\'ET\_10-1-2011-0471 and NIH
B330/479/11 ``Discrete and Convex Geometry''. 
The first and the last authors were partially supported by ERC
Advanced Research Grant no. 267165 (DISCONV), and the first author by
Hungarian National Research Grant K 83767, as well.
The second author was supported by the J\'anos Bolyai 
Research Scholarship of the Hungarian Academy of Sciences.
The third author was partially supported by MTM 2012-30719.
The fourth and fifth authors acknowledge partial support form CONACyT under
project 166306 and PAPIIT IN101912.

\end{document}